\documentstyle [11pt]{article}
\newcommand{\on}{\'{o}}

\def\bR{{\bf R}}   
\def\bC{{\bf C}}
\def\bZ{{\bf  Z}}      
\def\H{{\rm H}}
\def\bD{{\bf D}}        
\def\wylicz{
\list
 {\rm\arabic{enumi}.}{\settowidth\labelwidth{1.}%
\leftmargin\labelwidth
 \advance\leftmargin\labelsep
\parsep 2pt\itemsep 0pt\topsep 2pt
 \usecounter{enumi}}
 \def\newblock{\hskip .11em plus .33em minus .07em}
 \sloppy\clubpenalty4000\widowpenalty4000
 \sfcode`\.=1000\relax}

\def\szeroki{\xdef\stary{\normalbaselines}
\def\normalbaselines{\baselineskip 20pt
\lineskip 5pt \lineskiplimit 5pt}}
\def\wroc{\xdef\normalbaselines{\stary}}
\def\comm#1{\hbox{\hskip 5pt #1}}
\def\koniec{\hbox{\hskip 15pt $\Box$} \vskip 20pt}

\newcounter{pamietaj} \newcounter{temp} \newcounter{ch-temp}     
\newcounter{chp}
\def\tytul#1{\vskip 8pt
\setcounter{equation}{0}
\addtocounter{chp}{1}
\hbox to \hsize{\large\bf\arabic{chp}.\hskip 3pt #1\hfill}\vskip 10pt}
\def\u{\ifmmode {U(2)}\else $U(2)$\fi}
\def\su{\ifmmode SU(2)\else $SU(2)$\fi}
\def\so{\ifmmode SO(3)\else $SO(3)$\fi}
\def\sn{\hbox{\hskip 2pt\rm sn\hskip 2pt}}
\def\bz#1{\ifmmode \bZ_{#1} \else $\bZ_{#1}$\fi}

\def\M#1{{\cal M}_{#1}}
\def\df{D_f(\exp(t\sigma)\>) &=& D_f(}
\def\de{D_w(\exp(t\sigma)\>) &=& D_w(}
\def\kw#1{\sigma\cdot\sigma=-#1}

\let\duze=\sum 
\renewcommand{\sum}{\mathop{\displaystyle\duze}\limits}
\let\fff = \frac
\def\frac#1#2{\fff{\displaystyle#1}{\displaystyle#2}}

\newtheorem{theorem}{Theorem}[section]       
\newtheorem{lemma}[theorem]{Lemma}      
      
\newtheorem{cor}[theorem]{Corollary}

\newtheorem{definition}[theorem]{Definition}      

\newtheorem{nt}[theorem]{Note}

\newenvironment{note}{\begin{nt}\rm}{\end{nt}} 
\newenvironment{defn}{\begin{definition}\rm}{\end{definition}}

\title{\bf Immersed Spheres and Finite Type for \\
Donaldson Invariants \vskip.2in}
\author{\large Wojciech Wieczorek}

\addtocounter{section}{0}

\begin{document} 

\maketitle

\setcounter{equation}{0}
\section{Introduction}
\bigskip

Not that long time ago Donaldson invariants were the main tool in 
studying the four dimensional manifolds.
Even though these invariants have been overshadowed by the arrival of
the Seiberg-Witten invariants, they still are of interest to both
mathematicians and physicists.
For a simply 
connected smooth 4-manifold with $b_+ >1$ one can define a linear
function 
$$
\bD_w : {\bf A}(X) = {\rm Sym}^*\left( \H_0(X)\oplus \H_2(X)\right)
\to \bR
$$
for any $w\in\H^2(X,\bZ)$. By studying universal relations generated by
the presence of embedded surfaces Kronheimer and Mr\on wka \cite{k-str}
have shown that these invariants satisfy certain structure equations. To
make these relations even more compact, they have introduced a 
special class of manifolds.
\begin{defn}
 A manifold $X$ is of {\em $w$-simple type} for some $w$ if 
$D_w( (x^2 -4)
z) =0$ for any $z\in{\bf A}(X)$ and where $x$ is a generator of 
$\H_0(X)$.
\end{defn}

For manifolds with this property one can combine invariants of various
degree into a function 
 $\bD_w (\alpha) = D\left( (1+\fff{x}{2}) \exp(\alpha)\right)$.
This function $\bD$
has a very nice structure:
\begin{theorem}[Kronheimer, Mr\on wka \cite{k-str}]\label{str-th}
For a simply connected 4--manifold $X$ of simple type there are 
finitely many {\em basic} classes $K_1, \ldots , K_s\in\H^2
(X,\bZ)$ such that:
\begin{itemize}
\item[1.] { $K_i \equiv w_2 (X) \pmod 2$. }
\item[2.] {There are rational numbers $a_1,\ldots , a_s$ such that 
\begin{eqnarray}\label{wzor}
\bD_w = \exp(Q_X/2) \sum_{i=1}^s a_i \sinh {K_i}
&\quad&\mbox{
if $b_+\equiv 1 \pmod 4$, or}\\
\label{wzor2}
\bD_w = \exp(Q_X/2) \sum_{i=1}^s a_i \cosh {K_i} &\quad&
\mbox{
if $b_+\equiv 3 \pmod 4$.}
\end{eqnarray}
where $Q_X$ denotes the intersection form on $X$.
}
\end{itemize}
\end{theorem}

Soon after the discovery of this theorem,
Fintushel and Stern in \cite{fs} have shown that identical theorem 
can be
proven by studying embedded spheres. Moreover they proved that if a
 simply connected manifold is $w$-simple type for one $w\in\H^2(X)$
then it is of simple type for any other $w'$. 

It is still an open question whether there are simply connected
manifolds with $b_+>1$ that are not of simple type. 
The most important result that does
not use the assumption that a manifold is of simple type is
Fintushel-Stern blowup formula (\cite{blowup} and Theorem \ref{blup}),
which is a special kind of structure equation for manifolds containing
embedded spheres with self-intersection $-1$.
In an attempt to understand the structure equation in general case
 Kronheimer and
Mr\on wka  in non-published paper \cite{hyp} have defined a
manifold $X$ to be of {\em finite type} $r$ if it satisfies the condition
$D_w( (x^2-4)^rz) =0$ for some non-negative number $r$. The main
conjecture of their informal announcement is that all simply connected
manifolds with $b_+>1$ are of finite type. 
Recently Mu\~noz using Fukaya-Floer homology of $\Sigma\times S^1$
(Theorem 7.6 in \cite{munoz})
proved that when $X$ contains an embedded surface $\Sigma$ with genus $g$
and $\Sigma\cdot\Sigma=0$, then $X$ is of finite type with 
$$r=\left[\frac{2g+2}{4}\right]$$

In this paper instead of embedded surfaces we study immersed spheres.
First we focus on the relations between Donaldson invariants involving 
embedded spheres. We show how to write these relations  in a compact form. 
Next we move to studying immersed spheres. 
Our main result is the following:
\setcounter{pamietaj}{\value{equation}}
\begin{theorem}\label{moje}
\xdef\twierdzenie{
Let $\alpha$ be an immersed sphere with $p$ positive double points 
and any self-intersection $a$. Assume that there exist a
cohomology class $w$ such that $w\cdot \alpha = 1 \pmod 2$. Let $z$ 
be a product of classes perpendicular to $\alpha$.
     
For every $s=0,1,  \noexpand\ldots, p$ define the numbers:
$$
r = r(p,s) =\left[ \noexpand\frac{p+1-s}{2} \right]
$$
and
$$
k = k(a, s) = s - \left[\noexpand\fff{a+1}{2}\right] - 1
$$
$$
k_0 = \noexpand\cases{
 k & if $a$ is even\cr
k+1 & if $a$ is odd\cr}
$$

Then for every such $s$ there are the following structure equations:

$$
D_w\left( (x^2 - 4)^r \noexpand\cosh (t\alpha)z\right) = D_w\left( B^{-a} 
\noexpand\frac{(x^2 - 4)^r z}{(2-xq)^s} \cdot 
\sum_{i=0}^k  q^{i} Q'  \cdot c_i(\alpha)\right)
$$
and 
$$
D_w\left( (x^2 - 4)^r \noexpand\sinh (t\alpha)z\right) =  D_w\left( B^{-a}  
\noexpand\frac{(x^2 - 4)^r z}{(2-xq)^s} \cdot
\sum_{i=0}^{k_0} q^i Q \cdot d_i(\alpha)\right)
$$
where $c_i(\alpha)$ and $d_i(\alpha)$ are polynomials
of degree  $2i$ and respectively $2i+1$ on $\alpha$ 
(and some powers of $x$).}
\twierdzenie\ 
In the above $B$ and $S$ are the functions defined by Fintushel and 
Stern in their blowup formula \cite{blowup}, and we define 
$Q = (B/S)$ and $q= Q^2$.
\end{theorem}
As a consequence of this theorem we show that every simply connected 
 manifold containing an immersed sphere with $p$ positive double points 
and non-negative self intersection $a$  is of finite type with  
$$r=\left[\frac{2p+2-a}{4}\right]$$

We would like to thank Tom Mr\on wka and Eleny Ionel for helpful
discussions and encouragement.

\bigskip
\setcounter{equation}{0}
\section{The structure equations for embedded spheres}
\bigskip

The relations between various Donaldson invariants $D_w(\sigma^n)$
when a homology class $\sigma$ is represented by an embedded
sphere had been studied in \cite{fs} as well as in \cite{ww1}. 
In this
section we first review the basic definitions, after which we
state the main result for \u --invariants (Theorem \ref{gl-str}).
Since the proof of this Theorem is only a slight modification of our
previous result from \cite{ww1}, we include that proof in the 
Appendix. In this section, as a Corollary,
we will show how to find the coefficients of the general structure 
equations
in terms of some elliptic functions.

 To define the \u-- Donaldson invariants for
a  simply connected 4--manifold with $b_+\geq 3 $ odd 
we consider the compactified moduli space $\M {k,w} (X)$ of 
anti-self--dual (ASD)
connections on an \u--bundle $P$ over $X$ with $c_2
(P)=k\in\H^4(X,\bZ)$
 and first
Chern class $w\in\H^2(X,\bZ)$. Set $$
2d = {\rm dim } \M {k,w}(X) = 8k - 4 w^2 - 3(b_++1)
$$ 
Then there is a universal $SO(3)$ fibration ${\bf P}$ over 
$\M {k,w} (X)$ 
which gives rise to the homomorphism $\mu : {\rm  H}_i (X) \to 
{\rm H}^{4-i} (\M {k,w} (X) )$ given by 
$\mu (\sigma )= -\fff{1}{4} p_1 ({\bf P}) /  \sigma $. 
This allows one to define Donaldson invariants as linear maps     
$   
D_{d,w} : {\bf  A }^d(X)  ={\rm Sym}^d (\H_0 (X) \oplus \H_2 (X) )  
\to \bR    
$   
where the elements of $\H_i (X)$ have the degree $\fff{1}{2} 
( 4-i)$ and ${\bf  A }^d(X) $ is the set of elements of 
${\bf  A }(X)  = {\rm Sym}_* (\H_0 (X) \oplus \H_2 (X) ) $ 
having degree $d$. 
 The function $D_{d,w}$ assigns to the generator $x\in \H_0 (X)$ 
and the classes $\sigma_1, \ldots ,\sigma_r \in \H_2  (X)$ the 
number 
$$
D_{d,w} (\sigma_1 \cdots \sigma_r \cdot x^s ) = \langle \mu (\sigma_1)
 \ldots \mu (\sigma_r ) \cdot \mu (x)^s , [\M k (X) ] \rangle
$$
For simply connected manifolds with $b_+=1$ the same construction can
be performed, except that now the Donaldson invariant depends on 
metric on $X$.
The details of the above constructions can be found in \cite{dk}.

One combines the Donaldson invariants into a formal power series
$$
D_w (\exp(t\sigma)) = \sum D_w(\sigma^n)\frac{t^n}{n!}
$$

With the notation set above we can state the structure theorem for
embedded spheres:
\szeroki
\begin{theorem}\label{gl-str}
Let $X$ be a simply connected smooth 4-manifold that contains
an embedded 2-sphere $\sigma$ with self-intersection $\sigma
\cdot \sigma =s$. Let $z$ denote an arbitrary product of 
classes $z_i$ 
that are perpendicular to $\sigma$ and let 
$\varepsilon = w\cdot\sigma \pmod 2$
for some given cohomology class
$w$.  
Then there are universal functions $C_i = C_i(t,x;\varepsilon)$ 
such that when $s=-2k$:
$$
D_w(\exp(t\sigma)z\>) =\cases{ D_w\left((C_0 + C_1\sigma + \cdots + 
\widehat{\sigma^{2k-1}} +
C_{2k}\sigma^{2k})z\right)\hfill& 
when $\varepsilon = 0 $\cr
 D_w\left((C_0 + C_1\sigma + \cdots +
C_{2k-1}\sigma^{2k-1})z\right)\hfill& when $\varepsilon = 1 $\cr}
$$
and when $s = -(2k+1)$:
$$
 D_w\left((\exp(t\sigma)z\>)\right) = \cases{ D_w\left((C_0 + 
C_1\sigma + \cdots +
C_{2k}\sigma^{2k})z\right)\hfill & when $\varepsilon = 0 $\cr
 D_w\left((C_0 + C_1\sigma + \cdots +
\widehat{\sigma^{2k}} +
C_{2k}\sigma^{2k+1})z\right)& when $\varepsilon = 1$\cr}
$$
In the above the $\widehat{\sigma^{k}}$ means that the corresponding
 term does
not appear in the formula.
\end{theorem}
For  brevity we shall omit the class $z$ in  structure 
equations
that will follow, reminding about the condition $z\cdot\sigma = 0$
only when necessary. 

The first structure equation for embedded spheres that did not use the
simple type condition was Fintushel-Stern blowup formula. 
In  \cite{blowup} they  proved that when $\sigma$ is a class 
represented by an embedded sphere with self-intersection $-1$, then
for every $k$ there exist a function $B_k(x)$ such that:
\begin{equation}\label{blup}
D_w(\sigma^k)= D_w(B_k (x))
\end{equation}
for every $w\in\H^2(X)$ such that $w\cdot\sigma =0\pmod 2$.

The functions $B_k$ of  (\ref{blup}) when combined in the power 
series 
$B(x,t)= \sum_{k=0}^\infty \fff{t^k}{k!} B_k(x)$ satisfy:
\begin{equation}\label{bk}
B(x,t)= \exp\left(-\frac{t^2x}{6}\right) \sigma_3(t)
\end{equation}
where $\sigma_3$ is a particular quasiperiodic Weierstrass 
sigma-function  associated to the $\wp $--function $y$, which 
satisfies the differential equation
$$
(y')^2 = 4y^3 - g_2 y - g_3
$$
with $g_i$'s given by: 
$$
g_2 = 4\left(\frac{x^2}{3} -1\right)\hbox{, \hskip .1in} g_3= \frac{8x^3-36x}{27}
$$
(for details on elliptic functions see for example \cite{akh}).
Using the power series notation  we can rewrite (\ref{blup}) as 
\begin{equation} \label{bl-b}
D_w( \exp (t\sigma ))= D_w( B(t,x) )
\end{equation}

This formula is called the blowup formula (since after performing 
the
blowup on the manifold $X$, the exceptional divisor in $X\#
\overline{\bC
P^2}$ is represented by an embedded sphere with self-intersection 
$-1$). 
There is also a variant of the blowup formula for the twisted 
Donaldson
invariants, that originally was written as:
\begin{equation}\label{old}
D_{w+\sigma}( \exp (t\sigma ))= D_w( S(t,x) )
\end{equation}
where $S(t,x)= e^{-{t^2 x/6}} \sigma (t)$, and 
$\sigma (t)$ is
the standard Weierstrass sigma--function. 
To make our next formulas fit
into a single pattern we use the fact that $D_{w+\sigma}(\sigma)=
D_w$ and write (\ref{old}) as:
\begin{equation}\label{bl-s}
D_{w'}(\exp(t\sigma)) = D_{w'}(\sigma \cdot S(t,x))
\end{equation}
for $w'= w+\sigma$.

It turns out that we can express the functions $C_i(t,x,
\varepsilon)$ of Theorem \ref{gl-str} in terms of $B$ and $S$ defined 
above. For convenience we define
$$
\begin{array}{rcl}
\Delta&=& S'B-SB'\cr
Q&=& S/B\cr
q&=& Q^2
\end{array}
$$ 

With this notation we have:
\begin{theorem}\label{dokladnie}
Let $\sigma$ be a homology class represented by an embedded sphere
and set $n=-\sigma\cdot\sigma$. 
Like before, let $\varepsilon = w \cdot\sigma \pmod 2$.
Then modulo the kernel of $D_w(.)$ we have:
\renewcommand{\arraystretch}{2}
\begin{eqnarray*}
{\rm det}\> \left[
\begin{array}{c|c}
\cosh (t\sigma) & \left[ 1,\sigma^2,\ldots\right] \cr
\multispan{2}{\hrulefill}\cr
\left[ S^{2i} B^{n-2i}\right] & { W_e(n,\varepsilon)}\cr
\end{array}\right] &=&
 {\rm det}\> \left[
\begin{array}{c|c}
\sinh (t\sigma) & \left[\sigma,  \sigma^{3} ,\ldots\right] \cr
\multispan{2}{\hrulefill}\cr
\left[ S^{2i+1} \Delta B^{n-2i-3}\right] & { W_o(n,\varepsilon)}\cr
\end{array}\right] = 0 
\end{eqnarray*}
when $\varepsilon =0$, and:
\begin{eqnarray*}
{\rm det}\> \left[
\begin{array}{c|c}
\cosh (t\sigma) & \left[1,\sigma^2,\ldots\right]  \cr
\multispan{2}{\hrulefill}\cr
\left[ S^{2i} \Delta B^{n-2i-2}\right] & { W_e(n,\varepsilon)}\cr
\end{array}\right]& =&
{\rm det}\> \left[
\begin{array}{c|c}
\sinh (t\sigma) &  \left[\sigma,  \sigma^{3} ,\ldots\right]  \cr
\multispan{2}{\hrulefill}\cr
\left[ S^{2i+1} B^{n-2i-1}\right] & { W_o(n,\varepsilon)}\cr
\end{array}\right] = 0
\end{eqnarray*}
when $\varepsilon =1$.
The symbol $\left[ S^{2i+1} B^{n-2i-1}\right]$ denotes the 
column vector consisting of the entries indicated in brackets, where $i$
varies from $i=1$ to the maximum possible for which the power over $B$
is non-negative.
Similarly $ \left[ 1,\sigma^2,\ldots\right] $ denotes a row vector.
The $W_e(n,\varepsilon)$ and $W_o(n,\varepsilon)$
denote the Wronskian containing  even and correspondingly odd
 derivatives of functions
appearing 
in the first column of each matrix.
\end{theorem}
{\bf Proof:}
It turns out that all one needs to establish formulas like the ones
above is:
\begin{enumerate}
\item{The general structure equation as described in Theorem
\ref{gl-str}.}
\item{The blowup formulas (\ref{bl-b}) and (\ref{bl-s}).}
\end{enumerate}
We shall describe this procedure in the case
when $\sigma^2= -2k$ and $\varepsilon = 0$:

In this case 
there are  $2k$ functions $C_i(t,x)$ of Theorem \ref{gl-str}
to be found: the
$(k+1)$ coefficients with even powers of $\sigma$  and $(k-1)$ 
coefficients with
odd powers of $\sigma$. 
The structure equation is valid for any sphere of self-intersection
$-2k$ and any perpendicular class $z$, so we take any manifold $X$ 
with nontrivial Donaldson invariant and blow it up $2k$ times.
Let $e_1,
e_2, \ldots, e_{2k}$ denote the exceptional divisors, and let 
$\sigma =e_1+e_2+\cdots +e_{2k}$. Also let $w$ be a cohomology class
such that $w\cdot e_i = 0 $ for every $i$.

In order 
to find the first $(k+1)$ functions $C_i$ 
we  apply Theorem \ref{gl-str} for $w' = w+ e_1+e_2+\cdots +
e_{2i}$ where $i=0,1,\ldots , k$. Note that for such $w'$ and for
any odd number $s$,
$D_{w'}(\sigma^s)  = 0$. Thus we get $(k+1)$ equations of the
form:
$$
D_{w'}(\exp(t\sigma)) = 
D_w(S^{2i} B^{2k-2i}) =
D_{w'}\left((C_0 + C_2\sigma^2 + \cdots + 
C_{2k}\sigma^{2k})z\right)
$$
For each $r$ the number $D_w(\sigma^{2r})$ can be computed from the
blowup formula. In fact if we expand
$$
S^{2i} B^{2k-2i} = \sum_{j=0}^\infty m_j(x) \> t^j
$$
then $D_{w'}(\sigma^{2r}) = D_w( m_{2b}(x))$ which explains the 
appearance
of the Wronskian in Theorem \ref{dokladnie}. 

To find the remaining $(k-1)$ functions we consider
$\sigma$ as before and take  as an orthogonal class $z = e_1-e_2$.
Multiplying $D_w(\exp(te))$ by a class $e$ acts as differentiation, i.e.
$$
\matrix{
D_w(e\cdot\exp(te))&=& D_w(B'(t,x))\cr
D_{w+e}(e\cdot\exp(te))&=& D_w(S'(t,x))\cr}
$$ 
To ensure that we  get a nontrivial relation, the class $w'$ 
must
include precisely one of $\{ e_1, e_2\}$, e.g. $w' = w+ e_2+e_3+ 
\cdots +
e_{2i+1}$ for $i=1,2,\ldots , (k-1)$. 
Similarly like before we notice that for any even $n$
$$
D_{w'}(\sigma^n e_1) = D_{w'}(\sigma^n e_2) =0
$$ 
Then, plugging these in we get:
\begin{eqnarray*}
D_{w'}(\exp(t\sigma)(e_1 - e_2)) &=& D_w(S^{2i+1} B^{n-2i-3}
(B'S - S'B)) 
\hfill\cr
&=&
D_{w'}\left((C_1\sigma + C_3\sigma^3 + \cdots + 
C_{2k-3}\sigma^{2k-3})(e_1-e_2)\right)\cr
\end{eqnarray*}
We can arrange these equations in the 
Wronskian in a similar way as before. 

The proof of the other two cases is analogous.
\koniec

One can show that the functions $B$ and $S$ satisfy: 
$S(t,x) = t + O(t^3)$
and $B(t,x) = 1 + O(t^4)$. 
Thus for any $r$ and $n$:
\begin{equation}\label{s-power}
B^r(t,x) S^n(t,x) = t^n + O(t^{n+2})
\end{equation}
This shows that all Wronskians $W(n,\varepsilon)$ appearing in Theorem
\ref{dokladnie} are upper diagonal matrices with determinant equal to
one. Thus, for example, the first case of this Theorem may be 
rephrased as:
\begin{cor}
Let $X$ be a 4-manifold containing an embedded sphere $\sigma$ with
self-intersection $-2k$. Let $w\cdot\sigma =0$. 
Then there are polynomials $c_i(\sigma,x)$
of degree $2i$ in $\sigma$, such that:
$$
D_w(\cosh(t\sigma) ) = 
D_w\left( \sum_{i=0}^k c_i (\sigma ,x) S^{2i}B^{2k -2i}\right)
$$
\end{cor}
{\bf Proof:} Expand the corresponding matrix in Theorem \ref{dokladnie}
by the first column.
\koniec
 
Once we have the procedure of finding relations for embedded spheres,
we can ask a computer to find particular polynomials $c_i(\sigma)$
(we shall skip in our notation the dependence of these polynomials
on $x$.)
Next theorem lists several of such formulas, the last of them being 
of particular importance:

\begin{cor}\label{przyklady}
Let $\sigma$ be a homology class represented by embedded sphere and
$f$, $w \in\H^2(X,\bZ)$ such that $f\cdot\sigma = 1 \pmod 2$ and
$w\cdot\sigma = 0 \pmod 2$.

\noindent If $\kw 2$:
$$
\matrix{
\de B^2 + \sigma^2\frac{1}{2}S^2)\hfill\cr
\df \Delta + \sigma BS )\hfill\cr}
$$
If $\kw 3$:
\begin{eqnarray*}
\de B^3 + \sigma S\Delta + \sigma^2\frac{1}{2}BS^2)\hfill\cr
\df  B\Delta + \sigma S(B^2  +\frac{x}{3!}S^2)
+\sigma^3 \frac{1}{3!}S^3 \>)\hfill\cr
\end{eqnarray*}
If $\kw 4$
\begin{eqnarray}
\de (B^4 + \frac{1}{3} S^4) + \sigma SB\Delta + %
\sigma^2 \frac{S^2}{2}( B^2 +\frac{x}{3}S^2) 
+\sigma^4 \frac{1}{4!}S^4\>)\label{int4}\cr
\df \Delta \left(B^2+\frac{x}{2}S^2\right) +
\sigma SB(B^2 + \frac{x}{3!} S^2)+ \sigma^2 \frac{1}{2}S^2\Delta
\hfill + \sigma^3 \frac{1}{3!}S^3B \>)\hfill\nonumber\cr
\end{eqnarray}
\end{cor}

\begin{note} The formula for $-2$ sphere was 
first proved by R. Brussee and also has 
been known to other authors. The formula for $-3$ sphere and the 
class
$w$ was proven by the author in \cite{ww2}. At the time of writing 
that article we
could not handle formulas for spheres of self-intersection $-4$ 
or lower.
\end{note}

From the relation for $-4$ spheres 
 we get  the following formulas:
\begin{lemma}[\bf Double angle formulas]\label{double}
The functions $B$ and $S$ used in the blowup formula sa\-ti\-sfy:
$$
\matrix{D_w(B(2t))&=& D_w(B^4-S^4)\cr
D_w(S(2t))&=&2 D_w(\Delta BS)\cr
}
$$
Also we have:
$$
D_w(\Delta^2) = D_w(B^4 -xS^2B^2 +S^4)
$$
\end{lemma}
{\bf Proof:}
Let $\hat X$ denote $X$ blown up at one point and $e$ 
 an exceptional class of $\hat X$. For any natural number $n$ we 
have $D_X(x^n) = D_{\hat X} (x^n)$, thus it is sufficient to prove the 
Lemma for $\hat X$. The class $2e$ is
represented by embedded sphere with self-intersection $-4$. Note that
$f\cdot 2e \equiv 0 \pmod 2$ when $f=e$. Thus we can plug into the
first formula for $-4$ sphere $\sigma= 2e$ with $f\cdot \sigma =0$ to
get the ``double angle formula" for the function $B(2t)$, and then set
$f=e$ to get the formula for $S(2t)$. 
 In the last identity we blew up $X$ four times and 
used $\sigma= e_1 + \cdots + e_4$, $f=e_1+e_3$ and the perpendicular
class $z=(e_1-e_2)(e_3-e_4)$.
\koniec

\setcounter{equation}{0}
\section{Auxiliary lemmas}

In  computations involving the functions $B$ and $S$ it is crucial to
understand their quotient $Q = S/B$.
According to \cite{akh} we have:
\begin{equation}\label{q1}
Q(t) = \frac{\sigma (t)}{\sigma_3 (t)} = \frac{1}{\sqrt{e_1 - e_3}}
\sn (\sqrt{ e_1 - e_3} \cdot t, k)
\end{equation}
where:
\begin{eqnarray*}
e_1 = \fff{x}{6} + \fff{\sqrt{x^2-4}}{2} \hbox to .5in{}%
e_2 = \fff{x}{6} - \fff{\sqrt{x^2-4}}{2} \hbox to .5in{}%
e_3 = -\fff{x}{3}
\end{eqnarray*}
and:
$$
k^2 = \frac{e_2- e_3}{e_1-e_3} = \frac{x-\sqrt{x^2 -4}}
{x+\sqrt{x^2 -4}} =
\frac{\left( x-\sqrt{x^2 -4} \right)^2}{4} =
\frac{4}{\left( x + \sqrt{x^2 -4} \right)^2}
$$
With the above we can write (\ref{q1}) as: 
$$
Q(t) = \sqrt k \cdot \sn ( \frac{t}{\sqrt k}, k )
$$

The function $y=\sn (t)$ satisfies the differential equation:
$$
\left(\frac{dy}{dt} \right)^2 = (1-y^2)( 1- k^2 y^2)
$$
From this by simple computations we get:
\szeroki

$$
\matrix{ \left( Q'(t) \right)^2 &=& (1 - \frac{1}{k}\cdot Q^2)
(1- k\cdot Q^2)\cr
&=& 1 - x \cdot Q^2 + Q^4
}
$$
We could use this formula to prove the last statement in Lemma 
\ref{double}. In fact this shows that we can remove $D_w(.)$ from 
its statement.

We need one more statement about the functions $q^i$, which is immediate
from studying the Wronskian of these functions:
\begin{lemma}\label{niezalezny}
The functions $q^i$ are linearly independent.
\end{lemma}

For further computations we shall need some more elementary observations.

\begin{lemma}
Let $\sigma = \alpha + 2e_1 + 2e_2 + \cdots + 2e_k$ for some class 
$\alpha$ 
that is perpendicular to all exceptional divisors $e_j$. Let us fix 
a number 
$m\in \{0,1, \ldots, k\}$ and let $w $ be such 
that $w\cdot e_j = 0$ for every $j$.
Then for $w' = w+ \sum_{i=0}^m e_m$ we have:

\begin{equation}\label{sh-ch1}
D_{w'}(\cosh (t\alpha + \sum_{i=0}^k 2e_j)) =
\cases{D_w (\cosh (t\alpha) S(2t)^{i} B(2t)^{k-i}) & when $m$ is even. \cr
D_{w}(\sinh (t\alpha) S(2t)^{i} B(2t)^{k-i})& when $m$ is odd. \cr}
\end{equation}
and
\begin{equation}\label{sh-ch2}
D_{w'}(\sinh (t\alpha + \sum 2e_j)) =
\cases{D_{w}(\cosh (t\alpha) S(2t)^{i} B(2t)^{k-i})& when $m$ is odd.
 \cr
D_{w}(\sinh (t\alpha) S(2t)^{i} B(2t)^{k-i})& when $m$ is even. \cr}
\end{equation}

\end{lemma}
{\bf Proof:}
The result follows immediately from the relations:
$$
\matrix{
\cosh (\alpha + \beta ) &=& \cosh \alpha \cdot \cosh \beta +
\sinh \alpha \cdot \sinh \beta\cr
\sinh (\alpha + \beta ) &=& \cosh \alpha \cdot \sinh \beta +
\sinh \alpha \cdot \cosh \beta}
$$
combined with the blowup formulas (\ref{bl-b}), (\ref{old}).
\koniec

\begin{defn}
We shall call a polynomial $a_0 + a_1 x + \cdots + a_n x^n$ a {\em doubly 
monic of degree $n$} if both $a_0$ and $a_n$ are equal to $\pm 1$.
\end{defn}

\begin{lemma}\label{iloraz}
Let $f =f(q)$ be a doubly monic
polynomial of the degree $d$. 
Let $g$ be any polynomial of degree $k+d$. If 
$$
p = \frac{1}{f}\cdot g
$$
is a polynomial, then its degree does not exceed $k$.
\end{lemma}
{\bf Proof:}
By long division, we get
$$
g = f_1\cdot f + f_2
$$
Then a direct check 
of the coefficients of degree $k+d+1,$ $ k+d+2, \ldots$ in:
$$
p\cdot f = g = f_1\cdot f + f_2
$$
proves the result. Notice that $f$ needs to be monic for this
conclusion to hold, and its first coefficient needs to be $\pm 1$ 
for the formal inverse $1/f$ to make sense.
\koniec

\begin{lemma}\label{nww}
Let $f$ and $g$ be two doubly monic polynomials, and let $w_1$ and 
$w_2$ be some other two polynomials satisfying:
\begin{eqnarray}
f\cdot w_1 & = & g \cdot w_2\cr
f\cdot \phi_1 &+& g \cdot \phi_2 = C
\end{eqnarray}
for some polynomials $\phi_1$, $\phi_2$ and a constant $C$.
Then
$$
w_1 \cdot \phi_2 + w_2 \cdot \phi_1 = w_1 \cdot \frac{C}{g}
$$
\end{lemma}
{\bf Proof:} By direct computations:
$$
w_1 \cdot \phi_2 + w_2 \cdot \phi_1 = w_1\left( \phi_2 + 
(f/g)\cdot \phi_1\right) = \frac{w_1}{g} ( g \phi_2 + f\phi_1) =
C \frac{w_1}{g}
$$
\koniec
\wroc
\newcounter{axis}

\bigskip
\setcounter{equation}{0}
\section{The structure equation immersed spheres}
\bigskip

Now we are ready to prove our main result, which was stated in 
the introduction:
\hbox{}
\setcounter{ch-temp}{\value{section}}
\setcounter{section}{1}
\setcounter{temp}{\value{theorem}}
\setcounter{theorem}{\value{pamietaj}}
\begin{theorem}
\twierdzenie
\end{theorem}
\setcounter{theorem}{\value{temp}}
\setcounter{section}{\value{ch-temp}}

\setlength{\unitlength}{1mm}
\hbox to \hsize{\hfill\begin{minipage}{7cm}
\begin{picture}(60,60)
\put(5,10){\vector(1,0){60}}
\put(10,5){\vector(0,1){50}}

\multiput(10,10)(10,0){6}{\circle*{1}}
\multiput(20,20)(10,0){6}{\circle*{1}}
\multiput(30,30)(10,0){6}{\circle*{1}}
\multiput(40,40)(10,0){5}{\circle*{1}}
\multiput(18,8)(10,0){5}{\addtocounter{axis}{1}
\makebox(0,0)[tc]{\arabic{axis}}}

\multiput(10,10)(10,0){5}{\line(1,1){35}}
\put(60,10){\line(1,1){22}}

\put(70,10){\makebox(0,0)[cl]{$p$}}
\put(13,53){\makebox(0,0)[cl]{$s$}}
\put(45, 47){\makebox(0,0)[bc]{$r=0$}}
\put(60, 47){\makebox(0,0)[bc]{$\overbrace{\hbox to 1cm{}}^{%
\displaystyle r=1}$}}
\put(80, 47){\makebox(0,0)[bc]{$\overbrace{\hbox to 1cm{}}^{%
\displaystyle r=2}$}}

\end{picture}
\end{minipage}\hfill }
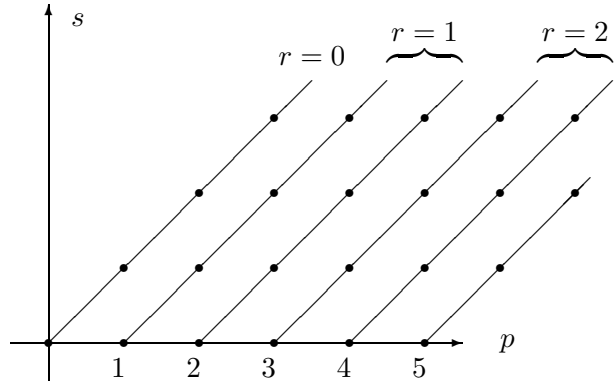
\centerline{Relation between $p$, $s$ and $r$ in Theorem \ref{moje}}
\vskip 10pt
\noindent {\bf Proof:}
We shall prove the Theorem by induction on $p$. 

For $p=0$ there is only one possible $s=0$, for which the conclusion of
the Theorem is the structure equation (\ref{gl-str}) with
 $r=0$. 

Thus we can assume that the Theorem is true for all $p' < p$. The
existence of the next $p+1$ structure equations we shall prove in the
following three steps:
\begin{enumerate}
\item{The $(p-1, s)$ - structure equation implies the $(p,s+1)$
structure equation.}
\item{The $(p-1, 0)$ - structure equation implies 
the  $(p,0)$ when $p$ is odd.}
\item{The $(p-2, 0)$ - structure equation implies 
the  $(p,0)$ when $p$ is even.}
\end{enumerate}

\vskip 5pt
\noindent{\bf Step 1.}
First let us fix $s\in\{0,1,\ldots, p-1\}$. 
Let $\beta = \alpha +2e$ be an immersed sphere with $p-1$ positive
double points obtained from $\alpha$ by blowing up one of its positive
double points. Notice that $r(p-1,s) = r(p,s+1)$, which we shall denote
simply as $r$. Denote $k = k(a,s+1)$.
Then for  $b = \beta\cdot \beta = a - 4$ we have:
$$
k(b,s) = s - b/2 -1 =  s - a/2 + 1 = k+1
$$
Similarly 
$$
 k_0(b,s) = k_0 + 1
$$
with $k_0 = k_0(a,s+1)$.
Let $w\in\H^2(X)$ be a cohomology class such that $w\cdot e = 0
\pmod 2$ (and of course $w\cdot\alpha = 1 \pmod 2$). Then
from the $(p-1, s)$ structure equation we have:
\szeroki
$$
\begin{array}{rcll}
D_w ((x^2 -4)^r \cosh (t\beta)) &=& D_w((x^2 -4)^r \cosh(t\alpha) 
B(2t)) & \hbox{by (\ref{sh-ch1}) and (\ref{bl-b})}\cr
&=& D_w\left(\frac{B^{4-a}}{(2-xq)^s} \sum_{i=0}^{k+1} c_i(\alpha)
\cdot q^i Q'\right)
& \hbox{by the inductive assumption}\cr
\end{array}
$$
which gives us:
\begin{equation}\label{step1-1}
D_w((x^2 -4)^r \cosh(t\alpha) (1-q^2)) = D_w\left(
\frac{B^{-a}}{(2-xq)^s} \sum_{i=0}^{k+1} c_i(\alpha)\cdot q^i Q'
\right)
\end{equation}
Note that $(w+e)\cdot(\alpha + 2e) = w\cdot\alpha \pmod 2$, thus 
similarly like before we obtain:
$$
\begin{array}{rrl}
D_{w+e}((x^2 -4)^r \sinh(t\beta)) = 
D_w((x^2 -4)^r \cosh(t\alpha) S(2t)) & \hbox{by (\ref{sh-ch2}) and 
(\ref{old})}\cr
 = D_w\left(\frac{B^{4-a}}{(2-xq)^s}
 \sum_{i=0}^{k_0+1} d_i(\alpha)\cdot q^i Q \right)
& \hbox{by the inductive assumption}\cr
\end{array}
$$
which after using $S(2t) = 2 QQ'$, dividing both sides by $Q$ and
finally multiplying by $Q'$ gives:
\begin{equation}\label{step1-2}
D_w((x^2 -4)^r \cosh(t\alpha)(1-xq +q^2)) = 
D_w\left(\frac{B^{-a}}{(2-xq)^s}
 \sum_{i=0}^{k_0+1} d_i(\alpha)\cdot q^i Q'\right)
\end{equation}
\begin{note}
Keep in mind that the functions $c_i(\alpha)$ and $d_i(\alpha)$ {\em
are not} quite the same as in the statement of the Theorem. In fact 
in this case one should use 
$$
D_w(\tilde{c_i}(\alpha)) = D_w(c_i ( \beta)) = 
D_w(c_i(\alpha + 2e) )
$$
In an attempt to simplify the notation we skip the extra ``tildes". 
\end{note}
\medskip

Now adding the  equations (\ref{step1-1}) and (\ref{step1-2}) 
(and setting $d_{k+2} = 0$ when $k_0 = k$) we  get the following relation, 
modulo the kernel of  $D_w(.)$:
\begin{equation}\label{step1-3}
(x^2 -4)^r \cosh(t\alpha) (2 -xq) = \frac{B^{-a}}
{(2-xq)^s}\left(
d_{k+2} q^{k+2} Q' +
\sum_{i=0}^{k+1}
(c_i (\alpha) + d_i(\alpha) )  \cdot q^i Q'  
\right)
\end{equation}
Dividing both sides by $(2 -xq)$ almost gives the required $(p,s+1)$
structure equation, except that the upper limit in the summation 
should be $k$ instead of $k+1$ for $a$ even, or $k+2$ when $a$ is
odd.

To eliminate extra terms, multiply (\ref{step1-1}) by $(1 -xq+q^2)$
and (\ref{step1-2}) by $(1-q^2)$ to get the relation:
\begin{equation}\label{step1-4}
D_w\left[\left(\sum_{i=0}^{k+1}c_i(\alpha)\cdot q^i\right)
(1 -x q + q^2)\right] = D_w\left[\left(\sum_{i=0}^{k+1}
d_i(\alpha)\cdot q^i + d_{k+2}q^{k+2}\right)(1-q^2)\right]
\end{equation}
Since the functions $q^i$ are linearly independent, the coefficient 
at $q^{k+4}$ must be equal, which means that $d_{k+2} =0$ also for 
$a$ odd. Comparing the coefficient at $q^{k+3}$ shows that
$c_{k+1} = - d_{k+1}$. This
implies that the $(k+1)$ coefficient in (\ref{step1-3}) is zero, 
thus proving this inductive step.

The proof for $\sinh(t\alpha)$ is similar (and even easier).

\vskip 5pt
\noindent{\bf Step 2.}
Let $r= r(p-1,0)$. If $p$ is odd, then:
$$
r(p,0) = \left[\frac{p + 1}{2}\right] = r + 1
$$
Like before, let $\beta = \alpha +2e$ be an immersed sphere obtained by
blowing up one of the positive double points of $\alpha$. Let us denote 
$k= k(a,0)$, so then 
$k(b,0) = k +2$. 

Similarly like before we get the relations:
\begin{equation}\label{step2-1}
D_w((x^2 -4)^r \cosh(t\alpha) (1-q^2)) = D_w\left(
\frac{B^{-a}}{(2-xq)^s} \sum_{i=0}^{k+2} c_i(\alpha)\cdot q^i Q'
\right)
\end{equation}
\begin{equation}\label{step2-2}
D_w((x^2 -4)^r \cosh(t\alpha)(1-xq +q^2)) = 
D_w\left(\frac{B^{-a}}{(2-xq)^s}
 \sum_{i=0}^{k_0+2} d_i(\alpha)\cdot q^i Q'\right)
\end{equation}
from which we get: 
\begin{equation}\label{step2-3}
D_w \left((1 -q^2)\sum_{i=0}^{k_0+2} d_i(\alpha)\cdot q^i \right) = 
D_w\left((1-qx + q^2)
\sum_{i=0}^{k+2} c_i(\alpha)\cdot q^i\right) 
\end{equation}
Direct check verifies the following identity:
\begin{equation}\label{magic1}
(1-q^2)\cdot (-2 - qx + x^2) + (1 - qx + q^2) \cdot ( -2 - qx) =
( x^2 -4 )
\end{equation}
Using this we get that $(x^2 -4)^{r+1} \cosh(t\alpha)$ 
modulo the kernel of $D_w(.)$ is equal to:
$$\begin{array}{ll}
(x^2 -4)^{r} \cosh(t\alpha) \left[ (1-q^2)\cdot (-2 - qx + x^2) + 
(1 - qx + q^2) \cdot ( -2 - qx)\right]& \hbox{by (\ref{magic1})}\cr
= B^{-a}\left[\left(\sum_{i=0}^{k+2} c_i(\alpha)
\cdot q^i\right) \cdot (2+ qx -x^2) + 
\left(\sum_{i=0}^{k_0+2} d_i(\alpha)\cdot q^i \right) \cdot ( 2+ qx)
\right]
\hfill&\hbox{by (\ref{step2-1}) and (\ref{step2-2})} \cr
=\left(\sum_{i=0}^{k+2} c_i(\alpha)\cdot q^i\right) \cdot
\frac{B^{-a}(x^2-4)}{(1-q^2)}\hfill&\cr
\end{array}
$$
In the last step we used  Lemma \ref{nww} for $f= (1-q^2)$ and 
$g=(1-qx+q^2)$ and relations (\ref{step2-2}) and (\ref{magic1}) 
as our assumptions.
Now applying Lemma \ref{iloraz} for the function $1 -q^2$
concludes the inductive step in this case. 

\vskip 5pt
\noindent{\bf Step 3.}
As before set $k=(a,0)$ and $r= r(p,0)$.
For $p$ even we have:
$$
r(p-2,0) = r-1
$$
The $p=0$ is our inductive assumption, thus we can assume that $p
\geq 2$. 
Let $\beta = \alpha + 2e_1 + 2e_2$ be an immersed sphere obtained from
$\alpha$ by blowing up two of its double points. Like above let $b=
\beta\cdot \beta$. Then
$$
k(b,0) = s-(a - 8)/2 -1 = k + 4
$$
From $(p-2,0)$ structure equation applied to $w$, $w+e_1$ and 
$w+e_1+e_2$ we get the following 
formulas:
\begin{eqnarray}
D_w((x^2 -4)^{r-1} \cosh(t\alpha) (1-q^2)^2) &=& D_w\left(B^{-a} 
\sum_{i=0}^{k+4}
a_i(\alpha)\cdot q^i Q'\right)\label{step3-1}\\
D_w((x^2 -4)^{r-1} \cosh(t\alpha) (1-xq+q^2)(1-q^2)) &=&  
D_w\left(B^{-a}\sum_{i=0}^{k_0+4}
b_i(\alpha)\cdot q^i Q'\right)\label{step3-2}\\
D_w((x^2 -4)^{r-1} \cosh (t\alpha)q (1-qx + q^2)) &=& 
D_w\left( B^{-a}\sum_{i=0}^{k+4}
c_i(\alpha)\cdot q^i Q'\right)\label{step3-3}
\end{eqnarray}
The relations (\ref{step3-2}) and (\ref{step3-3}) yield the following 
relation, modulo the kernel of $D_w(.)$:
$$
q\cdot\left( \sum_{i=0}^{k_0+4} b_i(\alpha)\cdot q^i \right) =
(1-q^2)\cdot\left( \sum_{i=0}^{k+4}
c_i(\alpha)\cdot q^i \right)
$$
By comparing the free term on both sides we obtain $c_0 = 0$, thus:
\begin{equation}\label{step3-4}
D_w((x^2 -4)^{r-1} \cosh (t\alpha) (1-qx + q^2)) = 
D_w\left( \sum_{i=0}^{k+3}
d_i(\alpha)\cdot q^i Q'\right)
\end{equation}
where $d_i(\alpha) = c_{i-1}(\alpha)$.
Putting together (\ref{step3-1}) and (\ref{step3-3}) we get that:
\begin{equation}\label{step3-5}
 \left(\sum_{i=0}^{k+4} a_i(\alpha)\right)(1-qx + q^2) =
\left(\sum_{i=0}^{k+3} d_i(\alpha)\right)(1-q^2)^2
\end{equation}
modulo the kernel of $D_w(.)$.
One can check directly that:
\begin{equation}\label{step3-6}
(1-q^2)^2 \cdot (x^2 -1 - qx) + (1 - qx +q^2) \cdot (-3 -2qx +q^2
+q^3 x) = (x^2-4)
\end{equation}
From which, similarly like in the previous step,
 we  conclude that $(x^2 -4)^{r} \cosh(t\alpha)$ is equal to:
$$
\begin{array}{l}
(x^2-4)^{r-1}\cosh(t\alpha)\left[
(1-q^2)^2 \cdot (x^2 -1 - qx) + (1 - qx +q^2) \cdot (-3 -2qx +q^2
+q^3 x)\right]\\
\hfill\comm{by (\ref{step3-6})}\\
=(x^2-4)^{r-1}\left[(x^2 -1 - qx)\cdot
\left(\sum_{i=0}^{k+4} a_i(\alpha)\right) 
+ (-3 -2qx +q^2 +q^3x)\cdot
\left(\sum_{i=0}^{k+3} d_i(\alpha)\right)\right]\hfill\\
\hfill\comm{by (\ref{step3-1}),(\ref{step3-4})}\\
 = \frac{(x^2-4)^r}{(1-q^2)^2}
\left(\sum_{i=0}^{k+4} a_i(\alpha)\right)\hfill\hbox{by 
Lemma \ref{nww}}\cr
\end{array}
$$
modulo the kernel of $D_w(.)$.
Now we can use Lemma \ref{iloraz} for the polynomial $(1-q^2)^2$ to 
conclude the induction. 
\koniec

\begin{theorem}
Let $X$ be a 4-dimensional manifold containing an immersed sphere
$\alpha$ with
$p$ positive double points and self-intersection $a\geq 0$. We 
also assume
that there exist a homology class $f$ such that $f\cdot\alpha = 1 
\pmod 2$. Then 
$$
D_w((x^2 - 4)^r) = 0
$$
for $r=\left[\frac{2p+2-a}{4}\right]$.
\end{theorem}
{\bf Proof:}
Assume first that $a$ is even. Then set
$$
s=\cases{a/2& if $a/2\leq p$\cr
p& otherwise\cr}
$$
With the above definition $r=\left[\frac{2p+2-a}{4}\right]$
and $k(a,s) \leq -1$. As a result the $(p,s)$ structure equation
of Theorem \ref{moje} has the form:
$$
D_w \left((x^2-4)^r \cosh(t\alpha)\right) = 0
$$
The coefficient at the free term $t^0$ gives the desired result.

When $a$ is odd, we first  blow up $\alpha$ at one of its {\em regular}
points, getting an immersed sphere with even self-intersection 
$b=a-1$. The result follows from the observation that for $a$ odd
$$
\left[\frac{2p+2-a}{4}\right] = \left[\frac{2p+2-(a-1)}{4}\right] 
$$
\koniec

\medskip
\appendix

\setcounter{equation}{0}
\section{Appendix: The proof of Theorem \ref{gl-str}}
\medskip

\subsection{Review of the gluing technique}
\bigskip

Let $X$ be an arbitrary smooth 4-manifold with $b_+>1$ that contains an
embedded sphere with self-intersection $-p$. In order to prove 
Theorem \ref{gl-str} we write the manifold $X$ as the sum $Y\mathop
\cup\limits_{L(p,1)} N$, where $N$ is the tubular neighborhood of the
sphere, and $L=L(p,1)$ is the lens space, equal to the boundary 
of $N$.

In $X$ we can identify a set isometric to $ (-\epsilon , 
\epsilon )\times L$ and take a sequence of metrics on $X$ that stretches
the length of $(-\epsilon , \epsilon )$ to infinity. 
In the limit we get two manifolds with cylindrical end
isometric to $\bR\times L$. We 
shall denote these manifolds by the same letters $Y$ and $N$. 
If one considers the anti-self-dual connections on \u\ bundles
over cylindrical end
manifolds, then the restrictions of these
connections to the bundles over $\{t\}\times L$ as $t\to \infty$ 
approach some
flat connection. More precisely, 
let 
$$\M {k,w}^x(Y) = \left(\{\hbox{ASD connections}\}\times \u\right)/{\cal G}
$$
 denote the space of based
 connections. Note that since the center $Z(\u)= S^1$ of \u\ acts 
trivially on the
space of connections, the fibration $\M {k,w}^x(Y) \to \M {k,w}(Y)$ 
is in fact an
\so -fibration.
It has been proved in \cite{mmrub} that there exists a  smooth 
boundary map 
$\partial ^o :\M{ k,w}^x (Y) 
\rightarrow {\cal R}(L)  $ from the based moduli space to the 
representation 
variety ${\cal R} (L)={\rm Hom} (\pi_1(L),U(2))/S^1$.
The representation variety ${\cal R} $
can be identified with gauge equivalence classes of
based flat connections. 

Let 
$$\chi = Map 
( \pi_1(L), G)/\sim$$
be {\em the
representation variety of $L$}.
 Here $\sim$ is the conjugation relation, 
i.e.
$\alpha\sim\beta$ if and only if there is some $g\in G$ such
that $\alpha= g\beta g^{-1}$. 
The  map $\partial^o$ descends to 
$\partial :\M {k,w} (Y) \rightarrow \chi (L)$. 

Fix $k$ and $w$ in $\M {k,w} (X)$. With the help of the boundary 
map $\partial^o$  we can define fibered product:
\begin{equation} \label{u-m}    
 U^{m}_{k_N,w } = \left(\M {{k_Y,w_Y}}^ x (Y,m) \times_  {{\cal R}} 
\M {{k_N,w_N}} ^x  (N,m)     \right) / SO(3)
\end{equation}
where $m\in\chi(L)$ and $(k_Y,w_Y)$, $(k_N,w_N)$ describe  Chern classes
of the restriction of the bundles to $Y$ and $N$ correspondingly. 
(So for example $k = k_N+k_Y$, thus since $k$ is fixed, we can drop $k_Y$ in
the notation above).
The space $\M {{k_Y,w_Y}}^ x (Y,m)$ denotes the space of those
connections on cylindrical end manifold $Y$ that are mapped onto 
$m\in\chi$ under the map $\partial$.
   
Let $X_l$ denote a manifold obtained from $X$ by stretching the cylinder
$ (-\epsilon , \epsilon )\times L$ to the length $l$.
Then we have the following theorem:

\begin{theorem}[see Theorem 4.3.1 of \cite{m-mrowka} or \cite{ww1}]
For $l$ large enough there exist a subset ${\cal Q}\subset 
\bZ\left[\fff{1}{p}\right]$, a subset $\chi_0\subset\chi(L)$ and a smooth map
$$
\mathop\coprod\limits_{(m, k_N)\in\chi_0\times {\cal Q}}  
U^{w} _{k_N,m}    
 \mathop{\rightarrow}\limits^\rho 
\M {k,w} (X_l)
$$ 
such that 
\begin{enumerate}
\item{$\rho$ is diffeomorphism onto its image.}
\item{If the term $\sigma^a x^b$ appears in one of formulas of Theorem
\ref{gl-str} and ${\rm dim }\M {k,w}(X) = 2a +4 b$, then the support of
the form $\mu^a(\sigma)\mu^b(x)$ is in the image of $\rho$.}
\end{enumerate}
\end{theorem}

\medskip
\subsection{Flat \u\  connections}
\bigskip

In \cite{ww1} we proved Theorem \ref{gl-str} for \su\ connections.
The main difference between \su\ connections   
and the \u\ connections is their representation variety. 

First let us refine the notion of the character variety. Since we 
study
connections on \u\ bundle with a fixed first Chern class $w$, then it
makes sense to define
$
{\chi}_w 
$
to be the set of equivalence classes of those flat connections that
``live" on the bundle with fixed first Chern class $w$. (In order
to simplify the notation we do not distinguish the class $w\in\H^2(X)$
and its restriction to $\H^2(L)$.) Similarly define ${\cal R}_w$.

Every matrix of \u\ can be
represented in the form:
$$
\left(\matrix{a&b\cr
-\bar b\phi & \bar a \phi\cr}\right)
$$
where $\vert a \vert^2 + \vert b \vert^2 =1$ and $\vert \phi 
\vert^2  =1$. 
Every element $g\in\u$ is conjugate (in fact by some \su\ matrix) to a
diagonal matrix 
$$\left(\matrix{\alpha&0\cr 0&\beta\cr}\right)
$$
where $\vert \alpha \vert^2= \vert \beta \vert ^2 =1$. 
Since $\pi_1(L) = \bz p = \H^2(L)$, we can identify $w$ with $\xi^a$ 
for some $a\in\bz p$ and $\xi$ a fixed primitive root of unity. 

Let $diag(k,l)$ denote the matrix
$\left(\matrix{\xi^k & 0 \cr 
0 & \xi^ l\cr} \right)$.
The relation:
$$
\left(\matrix{\xi^k & 0 \cr 
0 & \xi^ l\cr} \right) = 
\left(\matrix{0& -1 \cr
1 & 0\cr} \right)
\left(\matrix{\xi^l & 0 \cr 
0 & \xi^ k\cr} \right)
\left(\matrix{0& -1 \cr
1 & 0\cr} \right)^{-1}
$$
shows that the matrices $diag(k,l)$ and $diag(l,k)$ are conjugate.
Then:
$$
\chi_w = \{ diag(k,l) \>|\> k+l = a\}/\sim
$$
where $(k,l) \sim (l,k)$. There is a one-to-one correspondence between
elements of $\chi_w$ and the set
$$
\{i\in\bz{2p} \> |\> i = a \pmod 2 \}/ \pm
$$
In this correspondence the matrix $diag(k,l)\in\chi_w$ is send to 
$(k-l)\in\bz {2p}$. 

In the process of gluing it is also important to know what is the
inverse image of each component of $\chi_w$ in the fibration 
${\cal R}_w \to \chi_w$. This fiber $F_m$ is the stabilizer of 
$ m\in \chi_w$ divided by the central $S^1$.

When $\alpha \ne
\beta$ then the stabilizer is the subgroup of diagonal matrices, which
can be identified with $S^1\times S^1$. 
In case when $\alpha=\beta$ then the stabilizer is the whole group 
\u. Thus, after dividing by $S^1$ we have - like in the case of \su
\ connections - precisely two kinds of fibers: \so\ and $S^1$. In
\cite{ww1} we called the elements that have an \so\ as the fiber
to be the trivial elements of $\chi$. 

We can summarize this in the following lemma:

\begin{lemma}
Let $\chi_w$ denote the character variety of gauge equivalence classes of
flat \u\ connections on the bundle $E$ over $L=L(p,1)$ with first Chern
class $w$. We assume that the bundle $E$ extends to a tubular
neighborhood of embedded sphere $\sigma$ with self-intersection $p$.
 Then
if $p = -2k$, we have:
$$
\chi_{w} = \cases{ \{{\bf 0},2,\ldots ,{\bf 2k}\}&
if $w\cdot\sigma = 0 \pmod 2$\cr
 \{1,3,\ldots , 2k-1\}&if $w\cdot\sigma = 1 \pmod 2$\cr}
$$
and when $p = -(2k+1)$, then:
$$
\chi_{w}= \cases{ \{{\bf 0},2,\ldots , 2k\}&
if $w\cdot\sigma = 0 \pmod 2$\cr
\{1,3,\ldots , {\bf 2k+1}\}&if $w\cdot\sigma = 1 \pmod 2$\cr}
$$
In the above the bold face numbers indicate the trivial elements of the
character variety.
\end{lemma}     

\medskip
\subsection{Completion of the proof of Theorem \ref{gl-str}}
\bigskip

Now that we understand the character variety for \u\ connections, 
the proof of Theorem \ref{gl-str} is a slight modification of 
the corresponding structure 
theorem for \so\ 
connections (Theorem 1.3 \cite{ww1}). We shall
describe the main points of the proof, referring for the details
to \cite{ww1}.

Let $\M {k,w}(\bR\times L, [m,m'])$ denote the moduli space of ASD 
connections on the \u \ bundle over $\bR\times L$ with second Chern
 class $k$ and whose $\pm\infty$ limits are correspondingly 
$m,m'\in\chi_w$. 
Similarly let $\M {k,w} (N ,m)$ denote the moduli space of ASD
connections on cylindrical end manifold $N$, whose $\infty$ limits
are equal to $m\in \chi_w$.
Then we have the following:
\begin{theorem}
The dimension of $\M {k,w}(\bR\times L, [m,m'])$ is equal to 
\begin{equation}\label{dim-1}
dim \M {k,w}(\bR\times L, [m,m'])= 8k   + 2(m'-m) - 
\frac{2((m')^2- m^2)}{p} - s(m)
\end{equation}
where $s(m)$ is the dimension of the gluing parameter and is equal to:
$$s(m) = \cases{ 1 & if $m$ is non-trivial element of $\chi_w$\cr
3 &  otherwise\cr}
$$
The dimension of $\M {k,w} (N ,m)$ is:
\begin{equation}\label{dim-2}
dim \M {k,w} (N ,m) = 8k -3  + 2m - \frac{2m^2}{p} 
\end{equation}
\end{theorem}

\begin{theorem}\label{minimal}
For fixed $m$ and $m'$ let $k$ denote the minimal amount of energy 
for which the moduli space $\M {k,w}(\bR\times L, [m,m'])$ is 
non-empty. Then:
$$
k=\cases{\frac{(m')^2 - m^2}{4} & if $m' > m$\cr
\frac{(m')^2 - m^2}{4} + p(m-m') & otherwise \cr}
$$
In other words, the minimal dimension
for which the moduli space $\M {k,w}(\bR\times L, [m,m'])$ is 
non-empty is
$$
2(m' - m) - s(m)
$$
\end{theorem}

Let us assume  that $\sigma\cdot\sigma = -2a$ and
we want to prove that
there exist polynomials $C_i = C_i(x)$ such that:
$$
 D_w(\sigma^n\>) = \cases{D_w\left((C_0 + C_1\sigma + \cdots + 
\widehat{\sigma^{2a-1}} +
C_{2k}\sigma^{2a})z\right)\hfill& 
when $w\cdot\sigma = 0 \pmod 2$\cr
 D_w\left((C_0 + C_1\sigma + \cdots +
C_{2a-1}\sigma^{2a-1})z\right)\hfill& when $w\cdot\sigma = 1 
\pmod 2$\cr}
$$
Then similarly like in \cite{ww1} define the set
\begin{eqnarray}
{\cal J}_n&=&\{ (m,k)\> |\>  m\in\chi_w, \> k\in \bZ 
[\hbox{$\frac{1}{p}$}], \M {{k,w}} (N,m) \ne \emptyset, \cr
&&\hbox{\phantom{$\{ (m,k)\> |\>$} and }
 0< \>\dim  \M {{k,w}} (N,m) + s(m) \leq 2n \}   \nonumber 
\end{eqnarray}
For each $n$ this is a finite set, on which we 
define a partial order by saying that $ (m_1,k_1) \leq (m_2,k_2)$ 
if there is a nonempty moduli space $\M {{k_2 - k_1},w}
(\bR\times L , [m_1,m_2])$. Thus for example when $\sigma\cdot
\sigma = -6$ the set ${\cal J}_{10}$ is

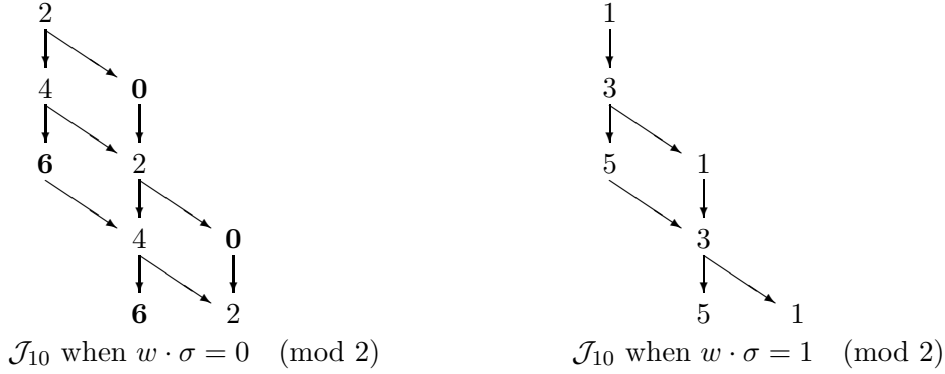
\setlength{\unitlength}{.5cm}
\begin{picture}(25,11.5)(-2,0)

\put(2,10){\makebox (0,0)[b]{2}} 
\put(2,8){\makebox (0,0)[b]{4}} 
\put(2,6){\makebox (0,0)[b]{\bf 6}} 
\put(4.5,8){\makebox (0,0)[b]{\bf 0}} 
\put(4.5,6){\makebox (0,0)[b]{2}} 
\put(4.5,4){\makebox (0,0)[b]{4}} 
\put(4.5,2){\makebox (0,0)[b]{\bf 6}} 
\put(7,4){\makebox (0,0)[b]{\bf 0}} 
\put(7,2){\makebox (0,0)[b]{2}}

\put(2,9.8){\vector(0,-1){1}} 
\put(2,7.8){\vector(0,-1){1}} 
\put(4.5,7.8){\vector(0,-1){1}} 
\put(4.5,5.8){\vector(0,-1){1}} 
\put(4.5,3.8){\vector(0,-1){1}} 
\put(7,3.8){\vector(0,-1){1}}  

\put(2,9.8){\vector(3,-2){1.9}} 
\put(2,7.8){\vector(3,-2){1.9}} 
\put(2,5.8){\vector(3,-2){1.9}} 
\put(4.5,5.8){\vector(3,-2){1.9}} 
\put(4.5,3.8){\vector(3,-2){1.9}}

\put(1,1){\makebox (0,0)[bl]{${\cal J}_{10}$ when 
$w\cdot\sigma = 0 \pmod 2$} }

\put(17,10){\makebox (0,0)[b]{1}} 
\put(17,8){\makebox (0,0)[b]{3}} 
\put(17,6){\makebox (0,0)[b]{5}} 
\put(19.5,6){\makebox (0,0)[b]{1}} 
\put(19.5,4){\makebox (0,0)[b]{3}} 
\put(19.5,2){\makebox (0,0)[b]{5}} 
\put(22,2){\makebox (0,0)[b]{1}}

\put(17,9.8){\vector(0,-1){1}} 
\put(17,7.8){\vector(0,-1){1}}  
\put(19.5,5.8){\vector(0,-1){1}} 
\put(19.5,3.8){\vector(0,-1){1}} 

\put(17,7.8){\vector(3,-2){1.9}} 
\put(17,5.8){\vector(3,-2){1.9}}  
\put(19.5,3.8){\vector(3,-2){1.9}} 

\put(16,1){\makebox (0,0)[bl]{${\cal J}_{10}$ when 
$w\cdot\sigma = 1 \pmod 2$} }

\end{picture}

In the above diagram we put only the $m$ coordinate of each
pair $(m,k)$. In order to obtain the  remaining $k$ coordinate 
we  
assign to each edge of the diagram the minimal energy defined in
Theorem
\ref{minimal}. Then $k$ coordinate of each vertex 
is the sum of those energies along any path from the
top of the diagram to that vertex.
Each vertex of $(k,m)\in{\cal J}$ represents one of the open 
set $U^w_ {k,m}$ defined in (\ref{u-m}), which cover the support of
$\mu(\sigma)^n$.
As we proved in \cite{ww1} the $\mu(\sigma)^n$ in each set $U^w_{k,m}$ 
can be computed in terms of  a polynomial $p_{n,m}(\mu(x),\mu(\sigma)
)$ whose degree in $\mu(\sigma)$ is at most $2m$.

For trivial a element $m$ of the character variety the dimension of 
the
gluing parameter $s(m)$ is by two bigger than $s(m)$ for all other 
elements. This is the reason for which we can delete the 
$\sigma^{2a-1}$ term in our relation.
\koniec

\end{document}